\documentclass[10pt,twocolumn]{article}

\usepackage{palatino}  
\usepackage[T1]{fontenc}

\usepackage{caption}
\usepackage{amssymb}
\usepackage{amsmath}
\usepackage{subfig}
\usepackage{graphicx}
\usepackage{wrapfig}
\usepackage{xspace}
\usepackage{textcomp}
\usepackage{sectsty}
\usepackage{sidecap}
\usepackage{eurosym}
\usepackage{pifont}
\usepackage{color}
\usepackage{mathpazo}
\usepackage{calc}
\usepackage{natbib}
\usepackage[american]{babel}
\usepackage{abstract}

\usepackage{geometry}
\definecolor{sectcol}{rgb}{0.63,0.16,0.16}
\definecolor{linkcol}{rgb}{0.459,0.071,0.294}

\usepackage[
    pdftex,
    pdftitle={A Tuned and Scalable Fast Multipole Method as a Preeminent Algorithm  for Exascale Systems},
    pdfauthor={Rio Yokota, L.~A. Barba},
    pdfpagemode={UseOutlines},
    bookmarks, bookmarksopen,bookmarksnumbered={True},
    colorlinks, linkcolor={linkcol},citecolor={linkcol},urlcolor={linkcol}
    ]{hyperref}

\makeatletter
\def\url@leostyle{%
  \@ifundefined{selectfont}{\def\UrlFont{\sf}}{\def\UrlFont{\small\ttfamily}}}
\makeatother
\def\rrr#1\\{\par
\medskip\hbox{\vbox{\parindent=2em\hsize=6.12in
\hangindent=4em\hangafter=1#1}}}

\newcommand{\nvidia}{n\textsc{vidia}\xspace}
\newcommand{\petfmm}{Pet\textsc{fmm}\xspace}
\newcommand{\petsc}{\textsc{pets}c\xspace}
\newcommand{\bigON}{$\mathcal{O}(N)$\xspace}
\newcommand{\bigOsq}{$\mathcal{O}(N^2)$\xspace}
\newcommand{\fmm}{\textsc{fmm}\xspace}
\newcommand{\pde}{\textsc{pde}}
\newcommand{\cpu}{\textsc{cpu}}
\newcommand{\gpu}{\textsc{gpu}}
\newcommand{\cuda}{\textsc{cuda}\xspace}
\newcommand{\simd}{\textsc{simd}\xspace}
\newcommand{\openmp}{\textsc{o}pen\textsc{mp}\xspace}
\newcommand{\mpi}{\textsc{mpi}\xspace}

\newcommand{\ME}{\textsc{me}}
\newcommand{\LE}{\textsc{le}}
\newcommand{\MM}{\textsc{m}2\textsc{m}\xspace}
\newcommand{\PM}{\textsc{p}2\textsc{m}\xspace}
\newcommand{\ML}{\textsc{m}2\textsc{l}\xspace}
\newcommand{\LL}{\textsc{l}2\textsc{l}\xspace}
\newcommand{\LP}{\textsc{l}2\textsc{p}\xspace}
\newcommand{\PP}{\textsc{p}2\textsc{p}\xspace}

\newlength{\up}
\newlength{\hup}
\sectionfont{\large\sf\bfseries\color{black}\vspace{0mm}}
\subsectionfont{\normalsize\it\bfseries\vspace{0mm}}
\subsubsectionfont{\normalsize\mdseries\itshape\vspace{-4mm}} 
\paragraphfont{\it}
\begin{document}

\pagenumbering{arabic}
\renewcommand{\thepage} { \arabic{page}}

\title{\sf \textbf {A Tuned and Scalable Fast Multipole Method as a Preeminent Algorithm  for Exascale Systems}}

\author{ \sf \textbf {Rio Yokota, L.~A. Barba}\\
\sf\normalsize  Mechanical Engineering Department, Boston University, Boston MA 02215}

\date{}

\twocolumn[
  \maketitle
  \begin{onecolabstract}
  
 Among the algorithms that are likely to play a major role in future exascale computing, the fast multipole method (\fmm) appears as a rising star. Our previous recent work showed scaling of an \fmm\ on \gpu\ clusters, with problem sizes in the order of billions of unknowns. That work led to an extremely parallel \fmm, scaling to thousands of \gpu s or tens of thousands of \cpu s. This paper reports on a a campaign of performance tuning and scalability studies using multi-core \cpu s, on the Kraken supercomputer. All kernels in the \fmm\ were parallelized using OpenMP, and a test using $10^{7}$ particles randomly distributed in a cube showed 78\% efficiency on 8 threads. Tuning of the particle-to-particle kernel using SIMD instructions resulted in $4\times$ speed-up of the overall algorithm on single-core tests with $10^{3}-10^{7}$ particles. Parallel scalability was studied in both strong and weak scaling. The strong scaling test used $10^{8}$  particles and resulted in 93\% parallel efficiency on 2048 processes for the non-SIMD code and 54\% for the SIMD-optimized code (which was still $2\times$ faster). The weak scaling test used $10^{6}$ particles per process, and resulted in 72\% efficiency on 32,768 processes, with the largest calculation taking about 40 seconds to evaluate more than 32 billion unknowns. This work builds up evidence for our view that \fmm\ is poised to play a leading role in exascale computing, and we end the paper with a discussion of the features that make it a  particularly favorable algorithm for the emerging heterogeneous and massively parallel architectural landscape.
    \vspace{3em}
  \end{onecolabstract}
]

\section{Introduction}

Achieving  computing at the exascale means accelerating today's applications by one thousand times. Clearly, this cannot be accomplished by hardware alone, at least not in the short time frame expected for reaching this performance milestone. Thus, a lively discussion has begun in the last two or three of years about programming models, software components and tools, and algorithms that will facilitate exascale computing. 

The hardware path to exascale systems, although not entirely certain, points to computing systems involving nodes with increasing number of cores, and consisting of compute engines with significant structural differences (for example, having different instruction set architectures), i.e.,  heterogeneous systems.
Among heterogeneous systems, those involving \gpu s are currently at the forefront and poised to continue gaining ground, until and unless a totally new and unexpected technology comes along. The latest Top500 list of the world's supercomputers (as of June 2011) has three \gpu-based systems ranked among the top five. Furthermore, on the Green500 list \cite[]{Green500} of the world's most energy efficient supercomputers (as of June 2011),  three out of the top five are also \gpu-based systems (not the same systems leading the Top500 list, however).
Even if a totally new technology were to come around the corner on our path to the exascale, the lessons learned from adapting our scientific algorithms to the challenges posed by \gpu s and many-core systems will better prepare us.  Features that are inevitable in emerging hardware models, as recognized in \cite{DarpaExascale2008}, are:  many (hundreds) of cores in a compute engine and programmer-visible parallelism, data-movement bottlenecks, and cheaper compute than memory transfers.  

Computational methods and software will have to evolve to function in this new ultra-parallel ecosystem. One of the greatest challenges for scalability (as well as programmability) is the growing imbalance between compute capacity and interconnect bandwidth.  This situation demands our investment in algorithmic research, co-developed with exascale applications so that the new architectures can indeed be exploited for scientific discovery at the highest performance.

On the other hand, it is clear that achieving exascale computing in the short time frame anticipated (the next decade) will require much more than hardware advances.  Developments in algorithms must play a leading role.  If history is to be our guide, algorithm developments have often provided great leaps in capability, comparable or higher than those offered by progress in hardware.  Thus, we can be fairly optimistic that, if we maintain a substantial research effort into both scientific algorithms and their implementation, the chances of success are higher.

It is interesting to note that many of the most successful algorithms of today are hierarchical in nature, such as the case of multi-grid methods.  Hierarchical algorithms often result in the ideal \bigON\ scaling, with $N$ representing the problem size.  A quintessential \bigON\ algorithm, which we focus on here, is the fast multipole method, \fmm.  This algorithm was introduced by \cite{GreengardRokhlin1987} for performing fast calculations of $N$-body type problems, such as are encountered in gravitational or electrostatic calculations.  These problems naturally have a computational complexity of \bigOsq for $N$ bodies (masses or charges), but the \fmm provides a solution in \bigON operations, to within a controllable accuracy. It is important to note that the \fmm\ can be used not only  as an $N$-body solver, as it was originally designed to be, it can also be used to solve any elliptic \pde\ \cite[e.g.,][]{EthridgeGreengard2001, ChengETal2006a,LangstonGreengardZorin2011}.

There are several groups actively contributing to the field of hierarchical $N$-body algorithms, and their high-performance implementation in both massively parallel \cpu\ systems and \gpu\ architecture.  These efforts have often been divided along two lines:  the $\mathcal{O}(N \log N)$ treecodes, with a strong following in the astrophysics community, and the \bigON\ fast multipole method, and its variants.
A notable player in the field is the kernel-independent version of the \fmm, a variant which does not make use of multipole expansions of the sources in a cluster but rather an equivalent density on a surface enclosing the cluster. The group working on this algorithm, known as \texttt{kifmm} \cite[]{YingBirosZorin2004}, has provided a continued stream of evidence of the high performance levels that can be attained with this family of algorithms. Most recently, they were awarded the coveted ACM Gordon Bell prize at the Supercomputing conference \cite[]{RahimianETal2010}, with a sustained performance of 0.7 Petaflop/s---this is merely $3.7\times$ less than the peak HPC Linpack benchmark of the same date, but achieved with a full application run (involving flow of deformable red blood cells). Thus, we are convinced that hierarchical algorithms, and the \fmm in particular, are especially poised to be leading players in the exascale computing future.

In our own recent efforts, we have achieved excellent scaling on \gpu\ clusters, with problem sizes that are in the order of billions of unknowns \cite[]{YokotaETal2011a,YokotaETal2011b}. This work has led to an extremely parallel \fmm, which is capable of scaling to thousands of \gpu s or tens of thousands of \cpu s. In the present study, we demonstrate the scalability of our new \fmm\ on a large \cpu-based system (Kraken, at the National Institute for Computational Science, University of Tennessee). We also perform single-node performance tuning by using \simd instructions and \openmp.   Our efforts are on a par with the most high-performing computations using \fmm, yet are fully independent. Thus, we contribute to the mounting evidence of the decisive role of \fmm and related algorithms in the quest for exascale. In view of this, we conclude our paper with a brief  discussion of the features of the \fmm\ that make it a  particularly favorable algorithm for the emerging heterogeneous, massively parallel architectural landscape. We hope that this will enrich the ongoing discussions about the challenges and opportunities to reach the milestone of exascale computing within a decade. 

\medskip

Given the magnitude of the challenge, we maintain that working in collaboration and under a model of openness should offer the best opportunities. Consistent with this view, all codes developed as part of this research effort are open source, and released at the time of publication. The full implementation of the \fmm used to produce the results of this paper, as well as some of our previous publications, is available via the website at \href{http://www-test.bu.edu/exafmm/}{http://www.bu.edu/exafmm/}.

\section{The fast multipole method, and its place in high-performance computing}

\subsection{Brief overview of the algorithm}

\begin{figure*}
	\centering
	{\includegraphics[width=0.75 \textwidth]{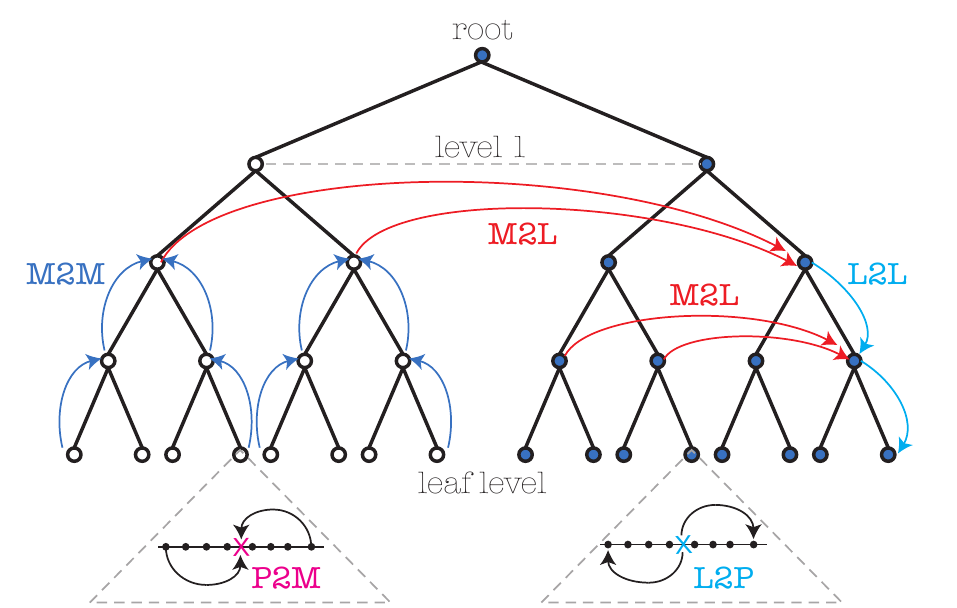}}
\vspace{-0.1cm}	\caption{\small Illustration of  \fmm algorithm components, with the  \emph{upward sweep} depicted on the left side of the tree, and the \emph{downward sweep} depicted on the right side of the tree.  The multipole expansions are created at the leaf level in the P2M operation, they are translated upwards to the center of the parent cells in the multipole-to-multipole (M2M) translation, then transformed to a local expansion in the M2L operation for the siblings at all levels deeper than level 1. The local expansions are translated downward to children cells in the L2L operation and finally, the local expansions are added at the leaf level and evaluated in the L2P operation.}
	\label{fig:tree5}
\end{figure*}

 Evaluating the interactions among $N$ sources of a force potential that has long-range effects in principle requires $\mathcal{O}(N^2)$ operations.  The canonical examples of this situation are gravitational and electrostatic problems. The fast multipole method (\fmm) of~\cite{GreengardRokhlin1987} is able to reduce the computational complexity of this problem to $\mathcal{O}(N)$ operations, and thus can have tremendous impact on the time to solution in applications.

\begin{figure*}
	\centering
	{\includegraphics[width=0.98 \textwidth]{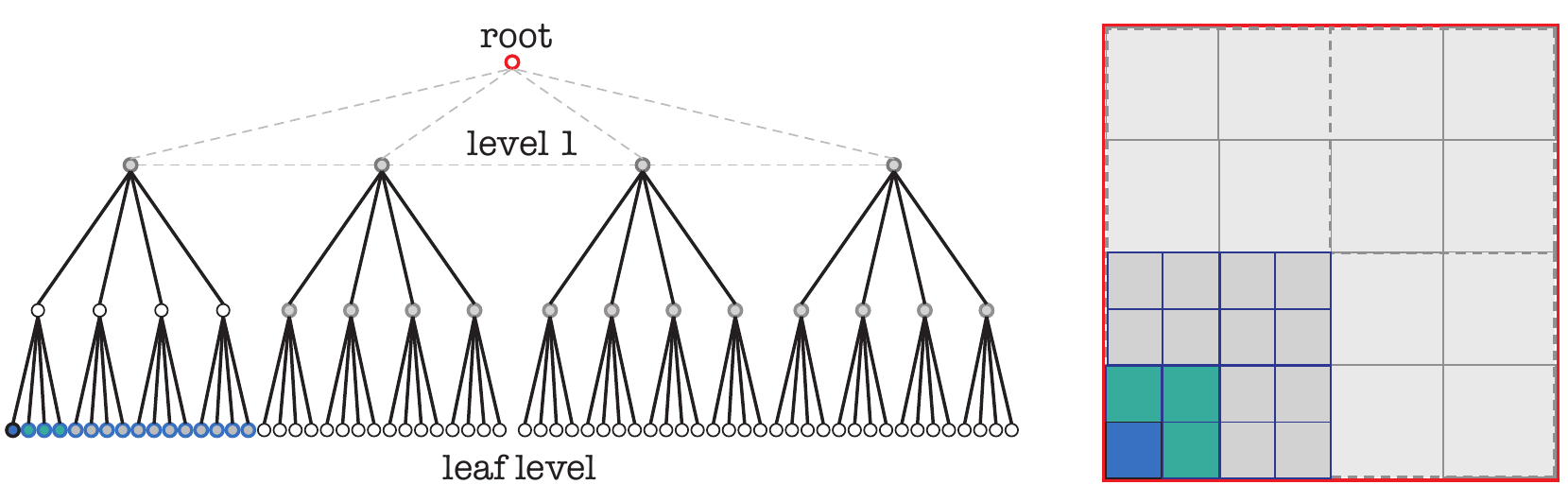}}
\caption{\small A three-level hierarchical decomposition of a 2D domain.  Tree nodes corresponding to domain boxes are shown in corresponding color. }
	\label{fig:MRdomain}
\end{figure*}

The \fmm  begins with a hierarchical subdivision of space, allowing a systematic categorization of spatial regions as either \emph{near} or \emph{far} from one another. 
 Different analytical approximations to the kernel that represents the interaction are used for pairs of points that are in near or far regions, correspondingly. 
These approximations are then transfered between different scales, up and down, using the subdivision hierarchy, eventually resulting in the assembly of a complete interaction force field at all of the points to a given approximation. 
A view of the complete algorithm is given in Figure~\ref{fig:tree5},  where the main computational steps are associated to a tree graph representing the hierarchical subdivision of a one-dimensional domain space.  The multi-resolution nature of the hierarchical space division, with its associated tree, is illustrated in Figure~\ref{fig:MRdomain} for the 2D case.

\begin{figure*}
\begin{center}
\includegraphics[width=0.95\textwidth]{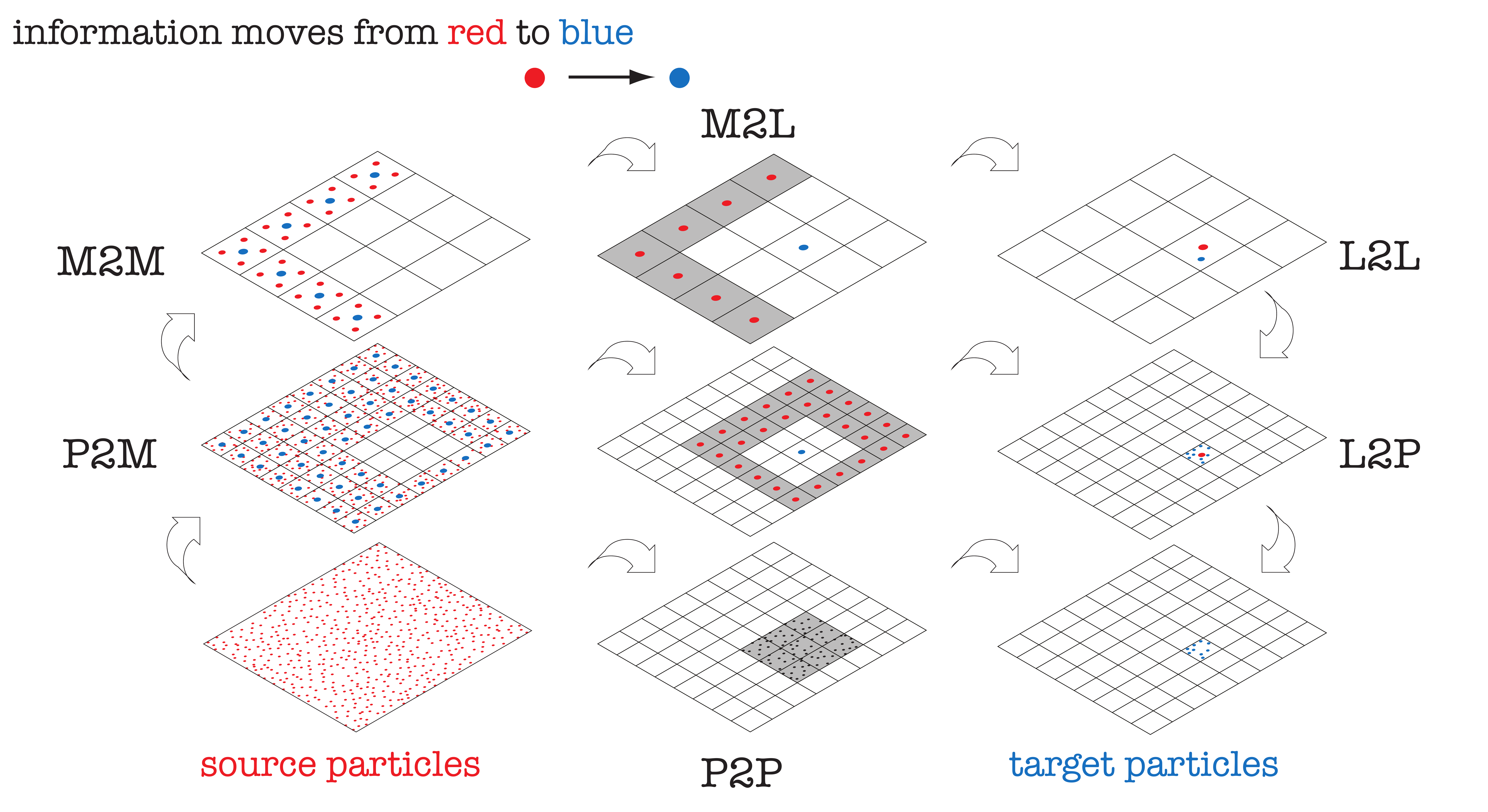}
\end{center}
\vspace{-1em}\caption{Flow of the \fmm calculation, showing all of the operations required in the algorithms: 
      \PM: transformation of points into \ME s (points-to-multipole);
      \MM: translation of \ME s (multipole-to-multipole);
      \ML: transformation of an \ME\ into an \LE\ (multipole-to-local);
      \LL: translation of an \LE\ (local-to-local);
      \LP: evaluation of  \LE s at point locations (local-to-point).}
\label{fig:flow_chart}
\end{figure*}

The \fmm\ algorithm defines the following mathematical tools. The Multipole Expansion (\ME) is a series expansion, truncated after $p$ terms, which represents the influence of a cluster of sources on far-away regions of space, and is valid at distances large with respect to the cluster radius. The Local Expansion (\LE) is also a series expansion of $p$ terms, valid only \emph{inside} a given subdomain, and used to efficiently evaluate the contribution of a group of \ME s on that subdomain.
In other words, the \ME s and \LE s are series (\emph{e.g}, Taylor series) that converge in different subdomains. The center of the series for an \ME\ is the center of the cluster of source particles, and it only converges outside the cluster of particles. In the case of an \LE, the series is centered near an evaluation point and converges locally. All the distinct operations that are required in the algorithm are illustrated in the flow diagram of Figure \ref{fig:flow_chart}.

The utilization of an aggregate representation for a cluster of particles, via the \ME, effectively permits a decoupling of the influence of the source particles from the evaluation points.  This is a key idea, resulting in the factorization of the computations of \ME s that are centered at the same point.  This factorization allows pre-computation of terms that can be re-used many times, increasing the efficiency of the overall computation. Similarly, the \LE\ is used to decouple the influence of an \ME\ from the evaluation points.  A group of \ME s can be factorized into a single \LE\ so that one single evaluation can be used at multiple points locally.  The combined effect of using \ME s and \LE s is a computational complexity of $\mathcal{O}(N)$ for a full evaluation of the force field.

\subsection{Our ongoing work on parallel \fmm}

Previously, our research group has led the development of an open-source parallel software, called \petfmm, which aims to offer a \petsc-like package for \fmm algorithms.  It is characterized by dynamic load-balancing based on re-casting the tree-based algorithm, illustrated in Figure \ref{fig:tree5}, into an undirected graph created using a computational model of the \fmm~\cite[]{CruzKnepleyBarba2010}. The spatial hierarchy is divided into many subdomains, more than the number of processes available, and an auxiliary optimization problem for the partition is constructed such that it minimizes both load imbalance and communication. The problem is modeled by a weighted graph, and solved using the partitioner ParMetis \cite[]{KarypisKumar98}.  
Strong scaling results showing the parallel performance of \petfmm were given in \cite{CruzKnepleyBarba2010}, with a test evaluating the field due to 10 million particles using up to 256 processors of a \cpu\ cluster.

\begin{figure}
\begin{center}
\includegraphics[width=0.50\textwidth]{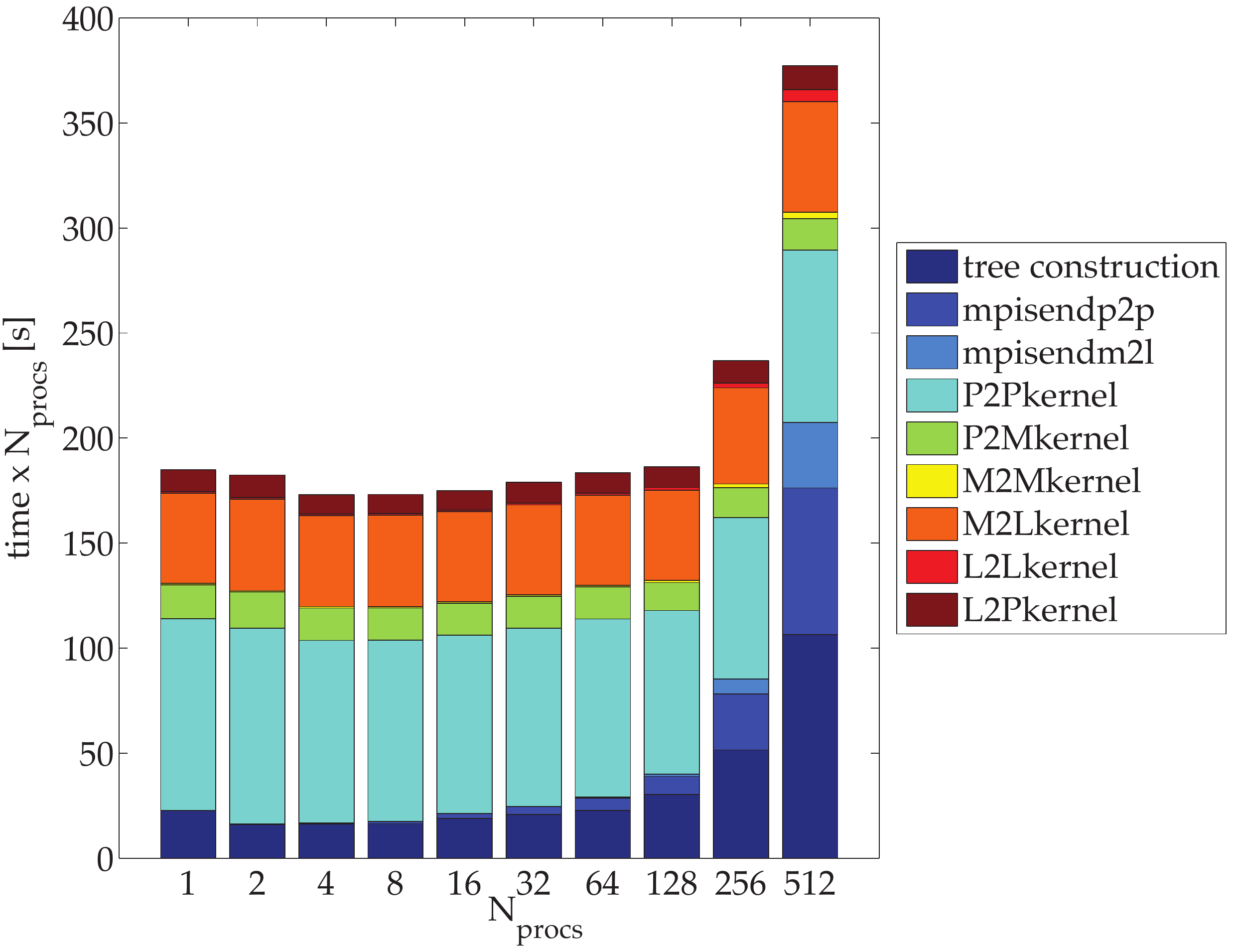}
\end{center}
\vspace{-0.5em}\caption{\textbf{Multi-\gpu\ strong scaling}, and timing breakdown of the different \fmm kernels, tree construction (done on the \cpu) and communications. Test problem: $10^{8}$ points placed at random in a cube; \fmm\ with  order $p=10$ and spherical harmonic rotation before each translation at $\mathcal{O}(p^{3})$ cost. Parallel efficiency is slightly above 1 for 4$-$32 processes, possibly due to cache advantage on the tree construction.  Parallel efficiency is 78\% on 256  and 48\% on 512 \gpu s. From \cite{YokotaETal2011b}.}
\label{fig:gpu}
\end{figure}

Recently, we have produced a multi-\gpu\ \fmm implementation, and performed new scalability studies.  
The parallel partitioning in this code is based on work-only load balancing, by equally distributing the Morton-indexed boxes at the leaf level; this is the standard and most commonly used technique for distributing $N$-body codes, first introduced by \cite{WarrenSalmon1993}. A host of detailed performance improvements were applied to accomplish good scalability on up to 512 \gpu s in our first study \cite[]{YokotaETal2011a}. We should note that achieving parallel scalability on a cluster of \gpu s is a greater challenge than scaling on a \cpu-only cluster. The high compute capability of \gpu s accelerates the computation and exposes the time required for communications and other overheads, which is exacerbated by the fact that, for multi-\gpu\ calculations, the \gpu s cannot communicate with each other directly. With current technology, it is necessary to send the data back to the host machine and communicate between nodes via MPI. 
Yet, our \textit{strong scaling} tests for the \fmm on multi-\gpu s demonstrated high parallel efficiency despite these challenges; see Figure \ref{fig:gpu}. On the Degima cluster\footnote{Degima is the \gpu\ cluster at the Nagasaki Advanced Computing Center, which at the time had two \nvidia \textsc{gtx295} cards per node (two \gpu s per card), and mixed \textsc{qdr/sdr} Infiniband network.}, parallel efficiency remains close to $1$ all the way up to 128 \gpu s, then drops to 78\% at 256, and to 48\% at 512 processes. The efficiency drop was ascribed to the configuration of the interconnect in the machine, and the number of \gpu s per node: up to 128 \gpu s, only one process was running on each node; for the case with 256 \gpu s, two processes were running on each node; and for the case with 512 \gpu s, four processes were running per node. This situation limits the effective bandwidth per process as bandwidth must be shared within the node. To alleviate this problem, we implemented a hierarchical all-to-all communication with \texttt{MPI\_Comm\_split}, splitting the MPI communicator into sub-sets of size $4$.  With $4$ processes per node, we are effectively splitting the communicator into an inter-node and an intra-node communicator. Thus, we first perform an intra-node \texttt{MPI\_Alltoallv}, then an inter-node \texttt{MPI\_Alltoallv}. We found that this could speed-up the communication nearly $4$ times at $512$ processes.  This is an example of applying a \emph{hierarchical parallel model} to effectively obtain performance from a heterogeneous system.

\subsection{High-performance \fmm computing}

Hierarchical $N$-body algorithms, including both the $\mathcal{O}(N\log N)$ treecode and the \bigON\ \fmm, have been visible players in the field of high-performance computing for many years. In fact, the top place in performance of the Gordon Bell award\footnote{\href{http://awards.acm.org/bell/}{http://awards.acm.org/bell/}} has been obtained several times with these algorithms. In 1992, the first place was awarded to \cite{WarrenSalmon1992} for a simulation of 9 million gravitating stars with a parallel treecode (sustaining 5 Gflop/s); some of the techniques used in that work remain important to this day. The same authors (with others) achieved the first place again in 1997 with a simulation of the dynamics of 322 million self-gravitating particles, at a sustained 430 Gflop/s \cite[]{WarrenSalmonETal1997}. That year, the award in the price/performance category also went to this team, who reported \$50/Mflop. 

More recently, we have seen treecodes and \fmm  figure prominently with Gordon Bell awards in the last two years. In 2009, the work by \cite{HamadaNarumiYokotaYasuokaNitadoriTaiji09} was recipient of the award in the price/performance category, achieving 124 Mflop/\$ (more than 6 thousand times ``cheaper'' than the 1997 winner) with a \gpu\ cluster. And last year, the first place in the performance category was achieved with an \fmm-based Stokes solver of blood flow which sustained 0.7 Pflop/s when solving for 90 billion unknowns \cite[]{RahimianETal2010}. To achieve this performance, a fast multipole method that already had undergone development to attain high parallel scaling \cite[]{LashukETal2009} was further optimized by means of explicit SSE vectorization and OpenMP multi-threading \cite[]{ChandramowlishwaranETal2010}. As a result, a full application code was obtained that was merely $3.7\times$ slower than the Linpack high-performance computing benchmark that is used to build the Top500 list\footnote{\href{http://www.netlib.org/benchmark/hpl/}{http://www.netlib.org/benchmark/hpl/}}.

\medskip

Given that the \fmm algorithm can achieve performance \emph{in applications} that is close to the Linpack benchmark, and that moreover it has excellent compatibility with massively parallel \gpu\ hardware, we estimate that it will feature prominently in the path to exascale. In the present contribution, we perform optimizations of our \fmm\ for single-core performance, and carry out scalability studies at large scale on \cpu-based systems. Added to our ongoing work on \gpu\ implementation, we keep the momentum going in the direction of exascale \fmm for the next decade.

\section{Intra-node performance optimization}

This section is dedicated to reporting a campaign of code optimization for single-node performance, including both \simd vectorization and multithreading within a node using \openmp. 
 All runs were performed on the Kraken Cray XT5 system of the National Institute for Computational Sciences (NICS) at the University of Tennessee, via TeraGrid access. 
 
 Each node of the Kraken system has two 2.6 GHz six-core AMD Opteron processors (Istanbul) and 16 GB of memory. The system is connected via a Cray SeaStar2+ router, and it has 9,408 compute nodes for a total of 112,896 compute cores.

\subsection{Multithreading with OpenMP}

\begin{figure}[t]
\begin{center}
\includegraphics[width=0.50\textwidth]{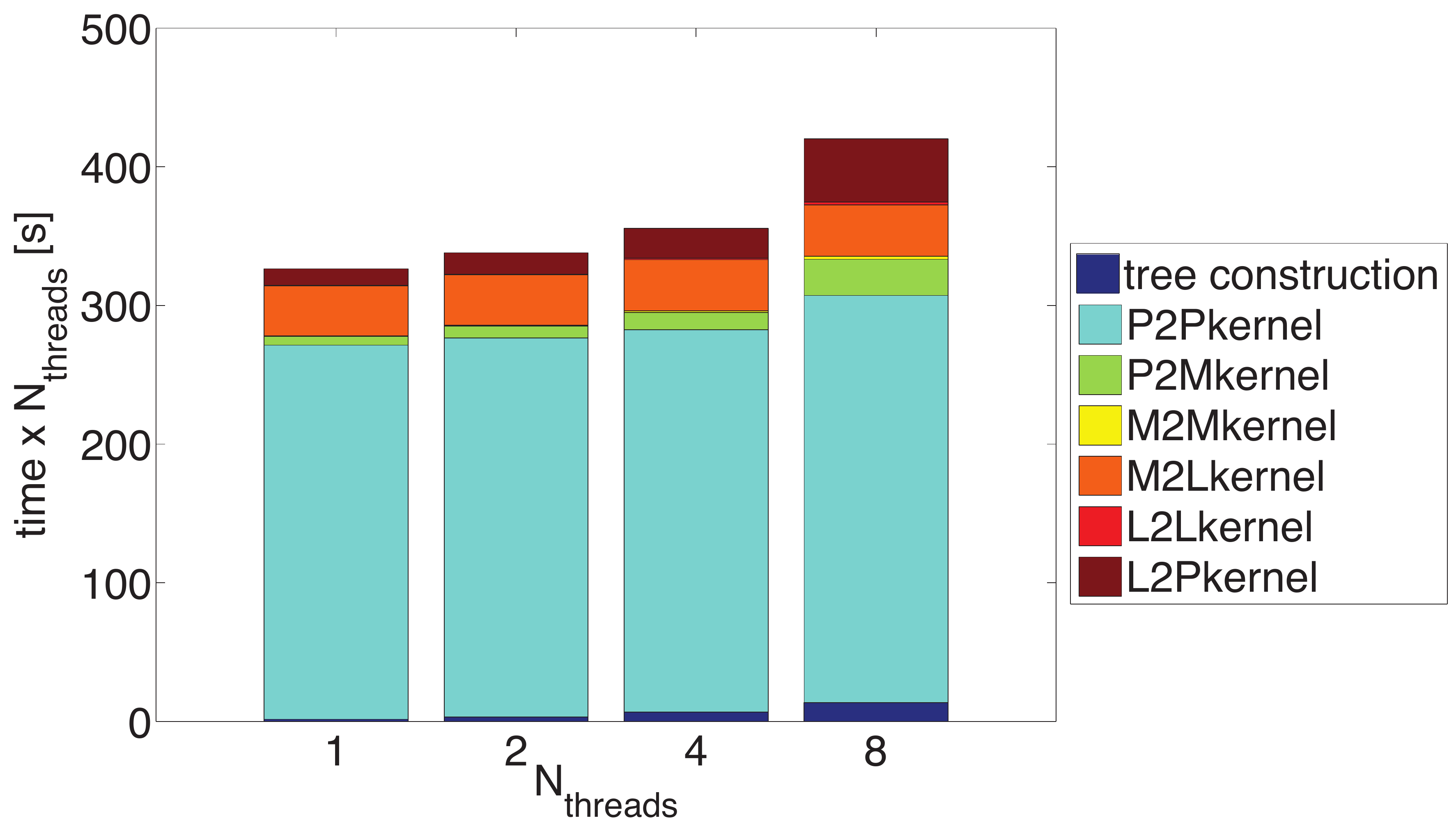}
\end{center}
\vspace{-0.5em}\caption{\textbf{\openmp strong scaling}, and timing breakdown of the different kernels, tree construction and communications. Calculation time is multiplied by the number of threads. Test problem: $N=10^{7}$ points placed at random in a cube; \fmm\ with  order $p=3$.  Parallel efficiency is 78\% on 8 threads.}
\label{fig:openmp}
\end{figure}

All kernels in the \fmm were parallelized using \openmp (the outermost loop for each kernel was explicitly parallelized by using the thread id). Results are presented in Figure \ref{fig:openmp}, consisting of the breakdown of the calculation time of the \fmm for different numbers of threads. The calculation time is multiplied by the number of threads for ease of comparison, i.e., bars of equal height would indicate perfect scaling. These tests were performed on a single node of Kraken using 1, 2, 4, and 8 cores, and the test problem was a set of $N=10^7$ particles randomly distributed in a cube.  In the legend on the right hand side, the \PP, \PM, \MM, \ML,\LL, \LP labels correspond to the individual stages of the \fmm as shown in Figure  \ref{fig:flow_chart}. The \PP kernel (the direct particle-particle interaction for the near field) takes up most of the calculation time, with the \ML (multipole-to-local translation between two distant cells) also taking a significant portion. These two kernels are balanced in an \fmm calculation by the number of levels in the tree, and the proportion of time taken by each will change as $N$ changes; in this case, with $N=10^7$, the situation shown in Figure  \ref{fig:openmp} offers the best timing and hence it is balanced.  The results in  Figure \ref{fig:openmp} show that the two main kernels scale quite well up to 8 threads. Some of the other kernels do not scale as well, which requires further investigation; however, this is not of primary importance since we do not need to scale to a large number of threads with the \openmp model at the moment.

The order of the \fmm expansions is set to $p=3$ for the tests in Figure \ref{fig:openmp}, which will give $3-4$ significant digits  of accuracy. We simply choose this value because it is common in astrophysics applications,  one of the main application areas for the \fmm, where high accuracy is not required. Using a larger $p$ would result in more work for the \PM, \MM, \LL, \LP, and especially the \ML kernel, which would require in turn using more particles per leaf cell to balance the calculation load of \PP and \ML kernels. This would affect the parallel scalability of the \PP kernel favorably, because of more data parallelism within each cell. Increasing $p$ would also increase the fine-grained parallelism of the \ML kernel because of more expansion coefficients to calculate, but would decrease the coarse-grained parallelism because the tree would become shallower for a given $N$.

\subsection{Vectorization with inline assembly}

\begin{figure}[t]
\begin{center}
\includegraphics[width=0.50\textwidth]{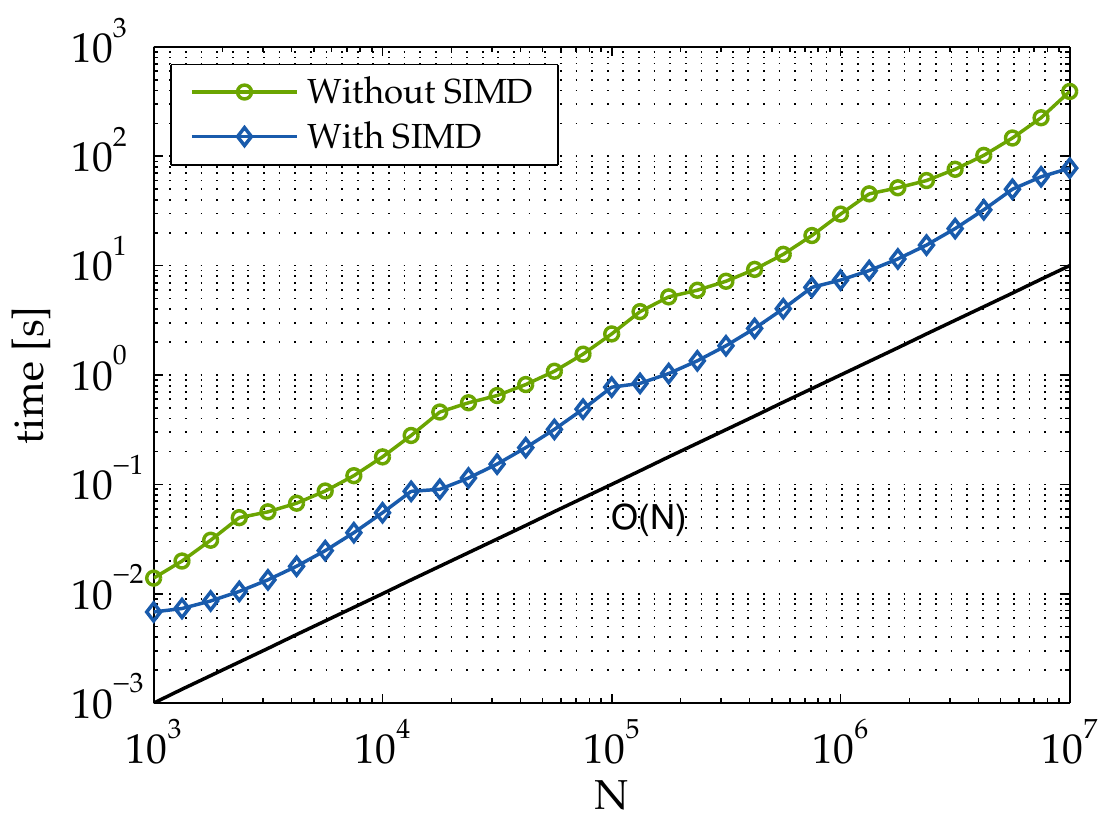}
\end{center}
\vspace{-0.5em}\caption{\textbf{Calculation time w/wo \simd} Test problem: $N=10^3-10^7$ points placed at random in a cube; \fmm\ with  order $p=3$. Only the \PP kernel is accelerated with \simd instructions. The overall speed up of the \fmm is approximately 4 times.}
\label{fig:simd}
\end{figure}

Recently, other researchers have dedicated substantial effort to single-node optimization and multi-core parallelization of an \fmm-type algorithm. \cite{ChandramowlishwaranETal2010} reported the first extensive study of its kind, which included not only SIMD vectorization and intra-node parallelization, but also various other optimizations of the \texttt{kifmm} algorithm. They reported a total performance improvement of $25\times$ for both optimizations and 8-core parallelization on Intel Nehalem \cpu s (plus differing performance improvements in other, less advanced, platforms). The speed-up from optimizations only on the overall algorithm is not specifically stated, but assuming perfect scaling on the 8 Nehalem cores,  we can infer a $3\times$ speed-up by the optimizations (including manual \simd vectorization) in double precision.

We have performed a similar optimization campaign using inline assembly and achieve a $4\times$ overall acceleration for a single-thread execution using single precision. Results are shown in Figure \ref{fig:simd}, consisting of the calculation time of the \fmm with and without \simd vectorization, for $N=10^3-10^7$ particles randomly distributed in a cube. At the moment, the optimization was performed for the \PP kernel only. The \PP kernel itself accelerates $16\times$ compared to an unoptimized double-precision kernel, but since the acceleration is limited to a certain part of the algorithm the overall acceleration is only $4\times$. We expect a further increase in performance when the other kernels are optimized using the same method. Nevertheless, even with the present implementation we are able to calculate $N=10^7$ particles in less than 100 seconds on a single core on Kraken (as shown in Figure \ref{fig:simd}).

It is worth noting that in citing the independent work on single-node optimization of the \texttt{kifmm} algorithm by \cite{ChandramowlishwaranETal2010}, we are not encouraging that our work be compared toe to toe in regards to performance. The \texttt{kifmm} algorithm is in many ways different from our \fmm, although it solves the same type of problems and shares a fundamental structure (e.g., hierarchical domain partition and associated tree structure). Quite the contrary, the point is that two completely independent efforts with hierarchical $N$-body algorithms are demonstrating the potential for applications based on these algorithms to benefit from excellent performance on multi-core systems.

\section{Parallel scalability on Kraken}

In this section, we report the results of parallel scalability tests on Kraken using the same \fmm code we described in the previous section.  Strong scaling is shown between 1 and 2048 processes and weak scaling up to 32,768 processes, with excellent parallel efficiency. For the strong scaling tests, we used 1 thread per process, and for the weak scaling tests we used 8 threads per process. The largest calculation is performed on a problem size of 32 billion particles, taking less than 40 seconds to complete.

The most recent results by researchers working with the \texttt{kifmm} algorithm and code base were reported in \cite{RahimianETal2010}, a work awarded the 2010 ACM Gordon Bell prize in the performance category. They report 84\% parallel efficiency for a strong scaling test between 48 and 3072 processes of the Jaguar supercomputer\footnote{Jaguar is a Cray XT5 system, very similar to Kraken.}. The study was carried out on an application of the \texttt{kifmm} algorithm to Stokes flow (red blood cells), and the reported strong scaling test used a problem of size $N=10^8$ unknowns. This is the same problem size that we used in our strong scaling test (below), but it must be noted that we use are using a different (Laplace) kernel. Moreover, the current literature points to many sophisticated features in the \texttt{kifmm} code (including: bottom-up parallel tree construction scheme, all-reduce-on-hypercube communication, FFT-based translation methods for the cell-cell interactions, among others), which are not utilized in our present software. Thus, it is not appropriate to make a performance comparison between these two very different variations of the \fmm and scaling studies. The point we would like to make here is that the \fmm is a highly ``strong-scalable'' algorithm, and therefore a truly suitable algorithm for the exascale era. The fact that two completely independent efforts both show similar excellent scalability is indeed promising.

\begin{figure}[t]
\begin{center}
\includegraphics[width=0.50\textwidth]{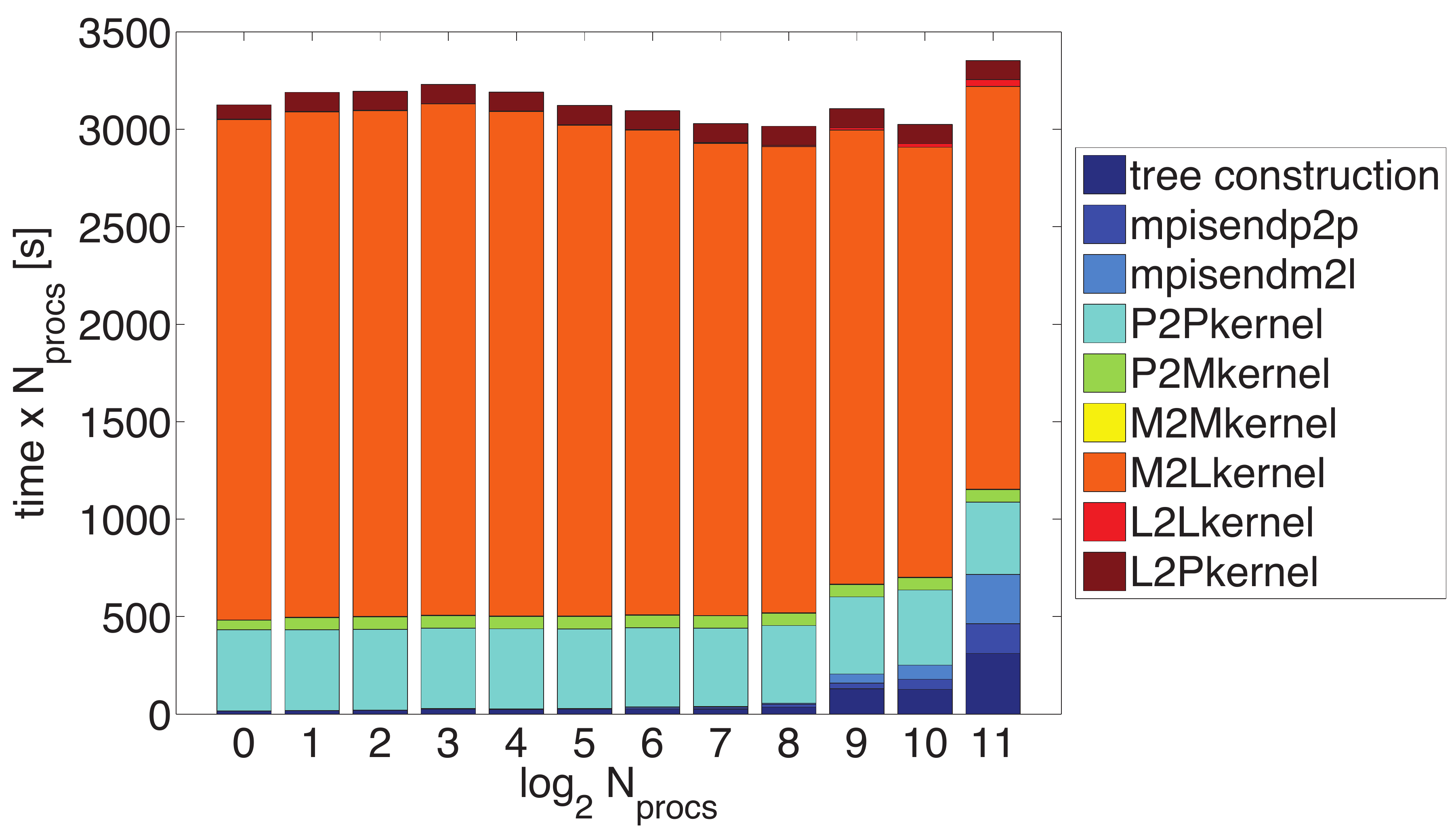}
\end{center}
\vspace{-0.5em}\caption{\textbf{\mpi strong scaling from 1 to 2048 processes}, and timing breakdown of the different kernels, tree construction and communications. Test problem: $N=10^{8}$ points placed at random in a cube; \fmm\ with order $p=3$. Calculation time is multiplied by the number of processes. Parallel efficiency is 93\% on 2048 processes.}
\label{fig:strong}
\end{figure}

We start by presenting our strong scaling results without  \simd vectorization. Figure \ref{fig:strong} shows the breakdown of each stage of the \fmm for a test consisting of $N=10^8$ points placed at random in a cube, with the order of multipole expansions set to $p=3$. The calculation time is multiplied by the number of \mpi processes ($N_{procs}$) for better visibility at large $N_{procs}$. The exponent of the number of processes is shown in the horizontal axis, e.g., the bar labeled with ``11'' represents $N_{procs}=2^{11}=2048$. The legend is similar to that of Figure \ref{fig:openmp} except we include the \mpi communication times  ``mpisendp2p" and ``mpisendm2l", which are the communication times for the \PP kernel and \ML kernel, respectively. All the other kernels do not require communication. The parallel efficiency at 2048 processes is in this case approximately $93\%$.

\begin{figure}[t]
\begin{center}
\includegraphics[width=0.50\textwidth]{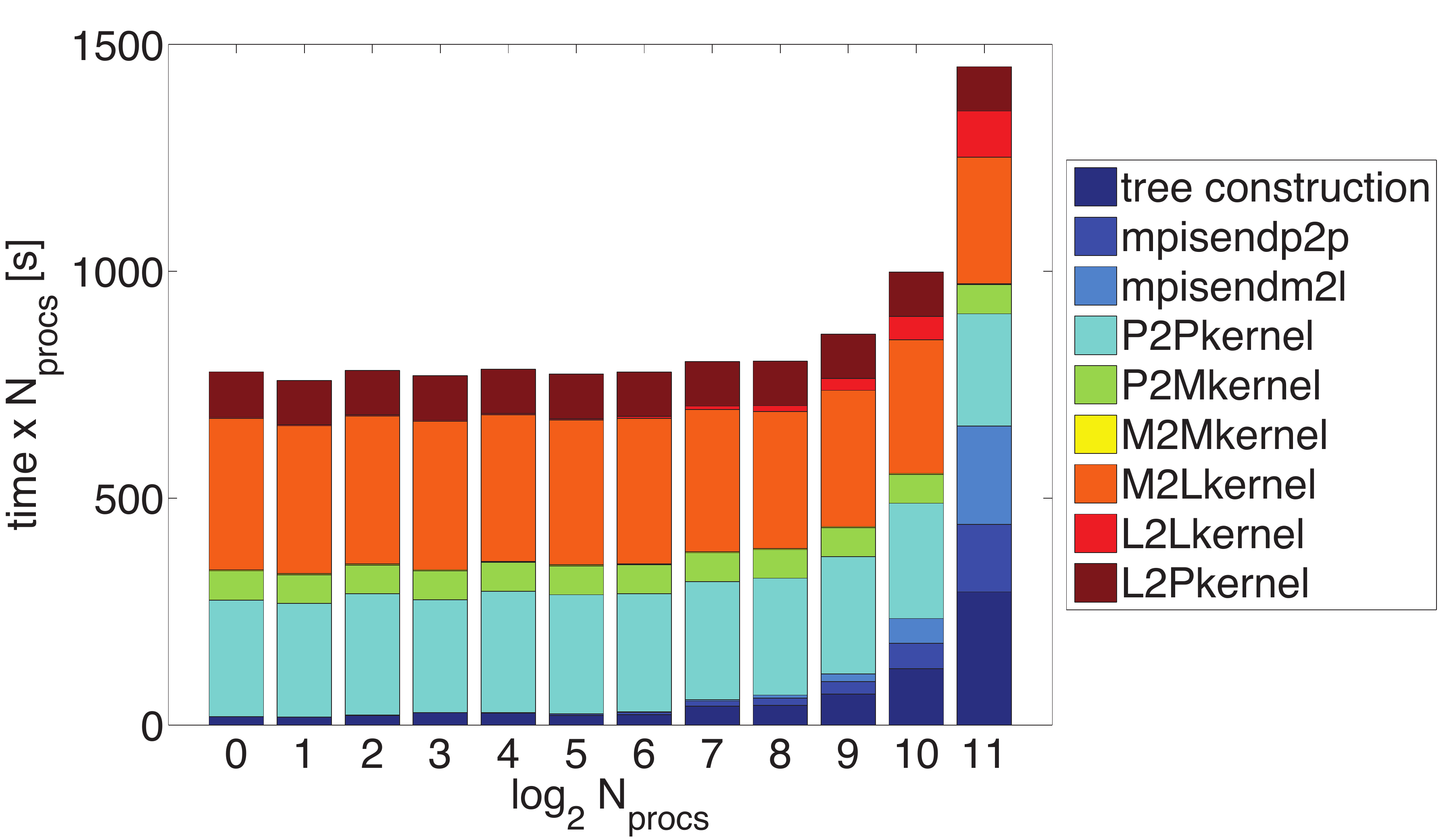}
\end{center}
\vspace{-0.5em}\caption{\textbf{\mpi strong scaling with \simd from 1 to 2048 processes}, and timing breakdown of the different kernels, tree construction and communications. Test problem: $N=10^{8}$ points placed at random in a cube; \fmm\ with order $p=3$. Calculation time is multiplied by the number of processes. Parallel efficiency is 54\% on 2048 processes.}
\label{fig:strong_sse}
\end{figure}

The results shown in Figure \ref{fig:strong} were obtained with the \fmm code \emph{without} single-node performance optimizations. The results for the same strong-scaling test using the \fmm \emph{with} the \simd kernels are shown in Figure \ref{fig:strong_sse}. As we have mentioned in the previous section, only the \PP kernel is optimized in our present code. Therefore, the tree is shallower and there are more particles in the leaf cell to keep the \PP and \ML balanced. Since the calculation time of the \PP kernel decreases significantly, the other parts start to dominate. The parallel efficiency is 54\% at 2048 processes, with about half of the total run time corresponding to the sum of tree construction time and time for communication between nodes and ---we must stress, however, that this is a strong-scaling test between 1 and 2048 processes (not offset scaling, as shown in some other works), and that communications are not overlapped.

By comparing Figures \ref{fig:strong} and \ref{fig:strong_sse}, we see that the calculation is approximately 4 times faster for the optimized version when the number of processes is small. However, since the optimized version does not scale as well, the speed-up is about 3 times at 1024 and only doubles at 2048 processes. These trends are similar to what we observed in the strong-scaling tests on the \gpu, where the acceleration of the \PP and \ML kernels causes the communication and tree construction to become the bottleneck.

\begin{figure}[t]
\begin{center}
\includegraphics[width=0.50\textwidth]{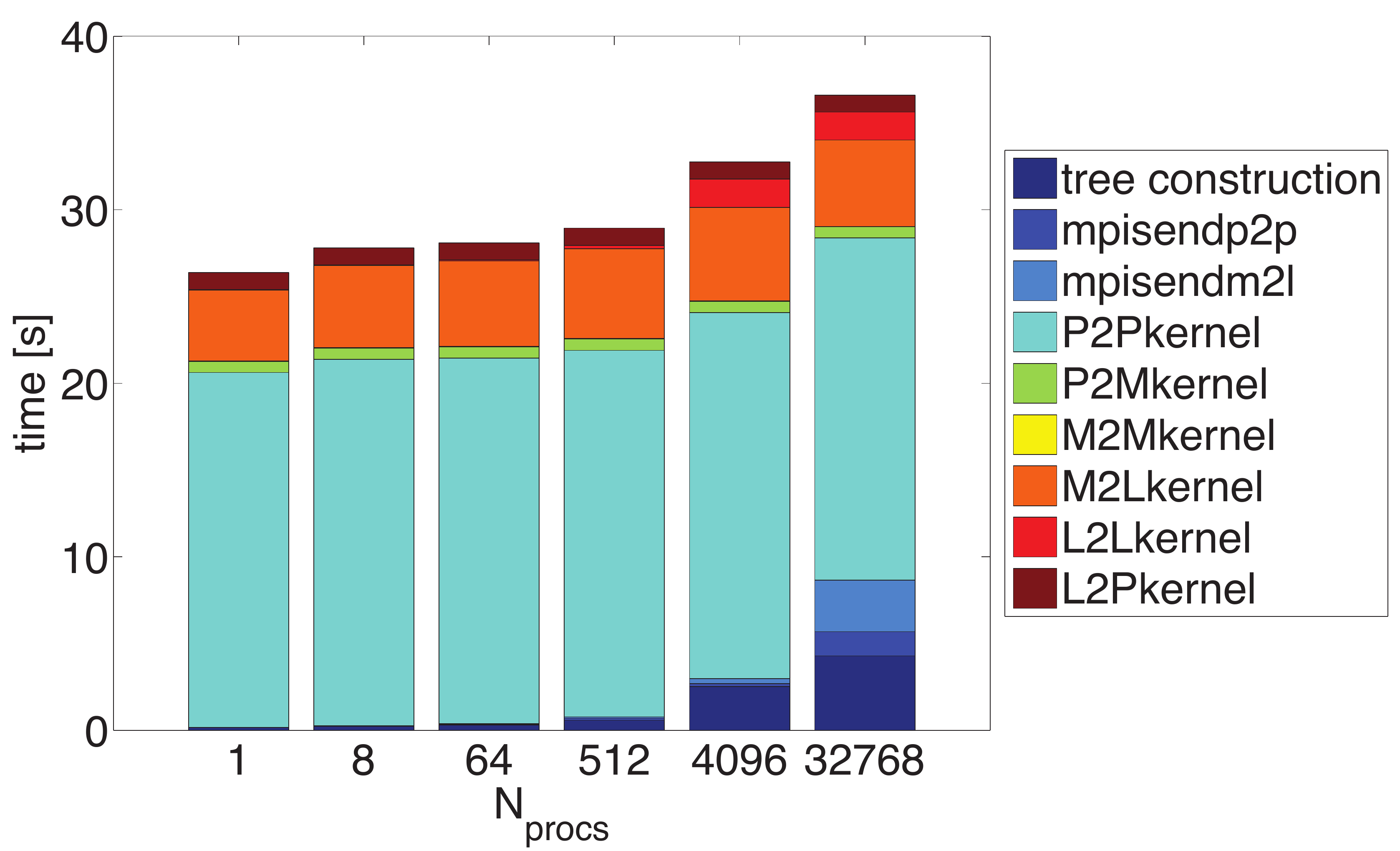}
\end{center}
\vspace{-0.5em}\caption{\textbf{\mpi weak scaling with \simd from 1 to 32768 processes}, and timing breakdown of the different kernels, tree construction and communications. Test problem: $N=10^{6}$ points per process placed at random in a cube; \fmm\ with order $p=3$. Parallel efficiency is 72\% on 32768 processes.}
\label{fig:weak}
\end{figure}

Figure \ref{fig:weak} shows weak scaling results on a test using $N=10^6$ particles per process and $p=3$. We achieve a parallel efficiency of 72\% on 32,768 processes. The largest problem size that was calculated consisted in \textbf{32 billion unknowns}, being solved in less than 40 seconds. It can be seen in Figure \ref{fig:weak} that the calculation times per process of both \PP and \ML kernels remain somewhat constant up to 32k processes.

\section{Suitability of the FMM for achieving exascale and future directions}\label{s:exafmm} 
 
The results shown here, added to our previous work and comparable recent work by other researchers, point to certain features of the \fmm that make it a particularly favorable algorithm for the emerging heterogenous, many-core architectural landscape.  In this section, we will make this view more clear, addressing specific exascale application characteristics identified in \cite{DarpaExascale2008}. The expectation for exascale is that in the order of a million terascale processors are required.
(Note that the results presented on Figure  \ref{fig:gpu} mean that we already have several million simultaneous threads running, if we consider the 32 non-divergent threads in a SIMD warp to be independent threads.) In recent feasibility studies, hypothetical exascale systems are proposed which would consist of 1024-core nodes, with each core clocking at 1 GHz, and a total of $2^{28}$ or $2^{30}$ such cores \cite[]{GahvariGropp2010,BhateleETal2011}. Analysis of representative algorithms can then be performed to determine the problem sizes required to reach 1 exaflop/s, and the constraints in terms of system communication capacity that would make this performance feasible. In  \cite{BhateleETal2011}, three classes of algorithms are considered by this analysis: pure short-range molecular dynamics, tree-based $N$-body cosmology, and unstructured-grid finite element solvers. \cite{GahvariGropp2010} similarly consider multigrid algorithms and fast Fourier transform. The conclusion of these studies is that the feasibility region for molecular dynamics (which is pleasantly parallel) and tree-based simulations is considerably less restricted than for other algorithms, having viable bandwidth requirements (in the order of $1-2$ GB/s). These analysis-based studies add to the experimental demonstrations of both our and other groups' ongoing work with \fmm to indicate great potential for this family of algorithms. Below, we explain some of the reasons for this potential.

\vspace{\up}
\paragraph{Load balancing.}Customarily, the \fmm is load-balanced via space-filling curves, either Morton or Hilbert.  The latter is an improvement because the length of each segment of the space-filling curve is more uniform, i.e., it never jumps from one space location to another far away.  These techniques balance only the computational work, without attempting to optimize communications.  We have demonstrated an approach  that optimizes both work and communications via a graph partitioning scheme in \cite{CruzKnepleyBarba2010}; however, this approach is not tested in large numbers of processors and extending it (by means of recursive graphs, e.g.) is an open research problem. In general, we believe that research into hierarchical load balancing techniques should be a high priority in our quest for exascale. The underlying tree-based structure of the \fmm points to an opportunity in this case.

\vspace{\up}
\paragraph{Spatial and temporal locality.} The \fmm is an algorithm that has intrinsic geometric locality, as the global $N$-body interactions are converted to a series of hierarchical local interactions associated with a tree data structure. This is illustrated in Figures \ref{fig:tree5} and \ref{fig:MRdomain}.  However, the access patterns could potentially be non-local.  An established technique is to work with sorted particle indices, which can then be accessed by a start/offset combination (this also saves storage because not every index needs to be saved).  Temporal locality is especially important on the \gpu, and is moreover programmer-managed.  We have implemented a technique that improves locality at a coarse level and is appropriate for \gpu s, namely, an efficient queuing of \gpu\ tasks before execution.  The queuing is performed on the \cpu, and buffers the input and output of data making memory access contiguous. On the other hand, locality at a fine level is accomplished by means of memory coalescing;  this is natural on the \fmm due to the index-sorting technique used. In conclusion, the \fmm is \emph{not} a ``locality-sensitive application'' \cite[]{DarpaExascale2008}. This is in contrast to stencil-type calculations, for example, where non-contiguous memory access is unavoidable.

\vspace{\up}
\paragraph{Global data communications and synchronization.}The overheads associated with large data transfers and global barrier synchronization are a significant impediment to scalability for many algorithms.  In the \fmm, the two most time-consuming operations are the \PP (particle-to-particle) and \ML (multipole-to-local) kernels.  The first is purely local, while the second effectively has what we can call \emph{hierarchical synchronization}.  The \ML operations happen simultaneously at every level of the hierarchical tree, without synchronization required \emph{between} levels.  This is in contrast, for example, to the situation in the multi-grid algorithm, which requires each level to finish before moving on to the next level. Among hierarchical algorithms, the \fmm appears to have especially favorable characteristics for ``exascaling''.

\subsection*{Future directions}

Based on the performance study reported in this paper, and also the features of the \fmm as mentioned above, we plan to continue improving our \fmm code and addressing some immediate tasks. First of all, strong scaling results in Figure \ref{fig:strong_sse} show that communication and tree construction become a bottleneck for continuing to achieve strong scaling beyond 2000 processes. Therefore, improvement of the tree construction phase is necessary, as well as overlapping of the local computation and communication. The tree construction can be accelerated by taking the particle distribution into account. For example, many applications in fluid dynamics and molecular dynamics have a uniform distribution of particles. For these cases, the tree structure is very easily formed by assigning a Morton index for a prescribed maximum level. It will be useful to have the capability to switch between uniform and adaptive tree construction, depending on the application. 
The efficient overlapping of communications and computations is the obvious next step. We are currently working on these improvements, and others, which will be reported in separate publications. 

These and other new features will be available in the \texttt{ExaFMM} code that we are releasing on the occasion of the Supercomputing conference in November 2011; for more information and access to the codes, visit \href{http://www-test.bu.edu/exafmm}{http://www.bu.edu/exafmm}.


\section{Conclusions}

The present contribution aims to build on the mounting evidence of the suitability of the fast multipole method, \fmm, and related algorithms, to reach exascale. The challenges to reach that milestone, which the community expects should be achieved by 2018, are varied and quite severe. They include serious technical hurdles, such as reducing the power consumption for computing by a significant factor, and dealing with frequent failures in systems involving millions of processors. While those challenges are recognized as the top concerns, the need for research into the appropriate algorithms and software for exascale is equally recognized as important. We maintain that \fmm is among the algorithms that will play a visible role in the exascale ecosystem.

We have conducted a campaign of single-node optimization and massively parallel scalability studies of our \fmm. Intra-node performance optimization is obtained both with multi-threading using \openmp, and with vectorization using inline assembly for the most costly kernel. \openmp parallelization of all kernels achieved 78\% efficiency in strong scaling with 8 cores, while inline assembly vectorization of only the \PP kernel provided a $4\times$ speed-up of the overall algorithm. 

Parallel scalability was studied  both in strong and weak scaling. Strong scaling between 1 and 2048 processes obtained 93\% parallel efficiency with the non-vectorized code, and 54\% efficiency with the SIMD-capable code. The SIMD instructions sped up the calculation by $4\times$ when using less than a thousand processes, and by $2\times$ when using 2048 proceses. In this last case, we start to see the inter-node communications and tree construction take a significant portion on the total runtime. However, even the vectorized code strong scales excellently up to one thousand processes, which is a situation rarely seen with other popular algorithms.
Weak scaling tests resulted in 72\% parallel efficiency on 32,768 processes for a test problem with a million points per process, and the largest calculation involved more than 32 billion points with a time to solution of under 40 seconds. 

Further work is ongoing to improve parallel efficiency in strong scaling by overlapping of communications and computations and enhancing tree construction. All such improvements are available on real time in the code repository, and the code is open for unrestricted use under the MIT License.

\section*{Acknowledgements}

This research was supported in part by the National Science Foundation (NSF) through TeraGrid resources provided by the National Institute for Computational Sciences (NICS) at University of Tennessee, under grant number TG-ASC100042. Additional support from NSF via grant OCI-0946441, from Office of Naval Research (ONR) via award \#N00014-11-1-0356, and from Boston University College of Engineering.  We also acknowledge guest access to the \gpu\ cluster at the Nagasaki Advanced Computing Center and to \textsc{Tsubame} 2.0 at the Tokyo Institute of Technology. LAB is also grateful of the support from n\textsc{vidia} via an Academic Partnership award.

\section*{Author Biographies}

\emph{Rio Yokota} obtained his PhD from Keio University, Japan, in 2009 and went on to work as a postdoctoral researcher with Prof.\ Lorena Barba at the University of Bristol and then  Boston University.   He has worked on the implementation of fast $N$-body algorithms on special-purpose machines such as \textsc{mdgrape}-3, and then on \gpu s after \cuda was released, and on vortex methods for fluids simulation.  In 2009, he co-authored a paper awarded the ACM Gordon Bell prize in the price/performance category, using \gpu s.  He joins the King Abdullah University of Science and Technology (KAUST) as a research scientist in September 2011.

\bigskip

\noindent\emph{Lorena Barba} is an Assistant Professor of Mechanical Engineering at Boston University since 2008.  She obtained her PhD in Aeronautics from the California Institute of Technology in 2004, and then joined Bristol University in the United Kingdom as a Lecturer in Applied Mathematics.  Her research interests include computational fluid dynamics, especially particle methods for fluid simulation; fundamental and applied aspects of fluid dynamics, especially flows dominated by vorticity dynamics; the fast multipole method and applications; and scientific computing on \gpu\ architecture. She was a recipient of the Engineering and Physical Sciences Research Council (EPSRC) First Grant in 2007, and was awarded an n\textsc{vidia} Academic Partnership grant in 2011.

\small
\bibliographystyle{jponew}

\end{document}